\documentclass[11pt]{amsart}
\usepackage{graphicx}
\usepackage{amssymb}
\usepackage{epstopdf}
\usepackage{amsmath}
\usepackage{amscd}
\usepackage{amsthm,amsfonts,amsxtra,amssymb,amscd}
\usepackage[all]{xy}
\usepackage{enumerate}

\usepackage{hyperref} 

\newtheorem{definition}{Definition}
\newtheorem{theorem}{Theorem}[section]

\newtheorem{lemma}[theorem]{Lemma}
\newtheorem{proposition}[theorem]{Proposition}
\newtheorem{corollary}[theorem]{Corollary}

\newtheorem{remark}[theorem]{Remark}
\newtheorem{remarks}[theorem]{Remarks}
\input xy

\xyoption{all}
\title{On the geometry of graph arrangements}

\author{C. De Concini}\address{Dip. Mat. Castelnuovo, Univ. di Roma La
Sapienza, Rome, Italy}\email{deconcin@mat.uniroma1.it}
\author{C. Procesi}\address{Dip. Mat. Castelnuovo, Univ. di Roma La
Sapienza, Rome, Italy}\email{procesi@mat.uniroma1.it}
\thanks{The authors are partially supported by the Cofin 40
\%, MIUR}
\begin{document}\begin{abstract} We use the results of \cite{dp}, \cite{dp1} to discuss the counting formulas of network flow polytopes and magic squares, i.e. the formula for the corresponding Ehrhart polynomial in terms of    residues. We also discuss a description of the big cells using the theory of non broken circuit bases.\end{abstract}
\maketitle
\section{Introduction} In this paper  we discuss  two topics which are complements of  the theory developed in \cite{dp}.   First an interesting example, the one of graph arrangements (in particular certain magic squares),  with applications to network flow polytopes.

Next (in an appendix)  we discuss a combinatorial approach to the determination of the big cells in a convex polyhedral cone, using the theory of non broken circuit bases, a concept introduced in  the theory of matroids and hyperplane arrangements.\medskip

Graph arrangements arise as follows: 
given a graph $\Gamma$ via its vertices $V$ and edges $L$  we 
fix a basis element $e_v$ for every vertex and, having chosen an orientation of the edges,  consider the set of vectors  $\Delta_\Gamma:=\{v_a:=e_{f(a)}-e_{i(a)}\}$
as $a\in L$ and $f(a), i(a)$ are the two vertices of $a$ (oriented).  

If  $V_\Gamma$  denotes the span of the vectors $v_a$, we have that 
the vectors $\Delta_\Gamma$ define  a hyperplane arrangement in   $V_\Gamma^*$.   This arrangement is clearly independent of the orientation chosen for the graph.\smallskip

The simplest example is the complete graph on $n+1$ elements which generates the arrangement of root hyperplanes for type $A_n$.\smallskip

Another special case is related to magic squares. 

Recall that a magic square of size $k $, is a square matrix of order $k$, filled with all the  integers $1,2,\dots,k^2$, with the property that the sums of the entries on each row or each column and the two diagonals is fixed.

In this paper we do not have really much to say about magic squares but rather we will study the weaker problem of counting  square matrices filled with integers with the property that the sums of the entries on each row or each column is a fixed number $n$.

This is  a special case of the theory of integral points in convex integral   polytopes among which are polytopes associated to oriented graphs (cf. \cite{BL}).

In our case  the  graphs given by choosing two disjoint sets $A,B$ with $m,n$ elements as vertices, and in which edges are all possible ones joining a vertex in $A$ with one in $B$.
\medskip

Recall that, given a  convex integral   polytope $\Pi$, the number of integral points in $n\Pi$, is a polynomial in $n$, called the  Ehrhart polynomial \cite{Er}, \cite{WV}.  We use the formulas developed in \cite{WV} as formulated in  \cite{dp} to compute this polynomial. The  formulas are based on a study of   an associated hyperplane arrangement,   this we explicit in this paper for graphs.

Let us recall the main points. In \cite{WV} the  Ehrhart polynomial  is given as follows. First a polytope is presented in the following way. We fix a set $\Delta:=\{\alpha_1,\dots,\alpha_N\}$ of vectors in a real vector space $V$. We assume that these vectors are all on the same side of some hyperplane. Then, given a vector $a$ in the positive cone spanned by the elements of $\Delta$, define \begin{equation}\Pi_a:=\{(a_1,\dots,a_N)\,|\, a_i\in\mathbb R^+,\ \sum_{k=1}^Na_i\alpha_i=a\}.\end{equation}
In order to count the integral points in $\Pi_a$ we pass to the corresponding hyperplane arrangement in  the dual space (we think of the $\alpha_k$ as linear equations defining hyperplanes in the dual space). We work then in the coordinate ring of the complement of the hyperplane arrangement.
 We make now a simplifying assumption (for the general theory   see \cite{SV1}).
 \begin{definition}
 A set $\Delta\subset\Lambda$ of integral vectors is said to be {\bf unimodular}, if all the subgroups spanned by subsets of $\Delta$  are direct summands of $\Lambda$.
 \end{definition}
  The final formula that we use is the one developed in \cite{dp} and reproduced here in \S   3 (5), (6).

The first half of this paper is just a reminder of the notations and results of \cite{dp},\cite{dp1}. The reader familiar with this work can pass directly to the last section where we illustrate the general theory in the case of graph arrangements.

 \subsection{Notations} With the notations of the introduction,  let $U$
be a complex vector space of dimension $r$,    $\Delta\subset U$   a
totally ordered finite set of vectors
$\Delta=\{\alpha_1,\ldots ,\alpha_m\}$. These vectors are the linear
equations of a hyperplane arrangement in $U^*$. We also
assume that $\Delta$ spans $U$ and any two distinct elements in $\Delta$
are linearly independent.

 From these data one constructs
the partially ordered set of subspaces obtained by intersection of the given
hyperplanes and the open set $\mathcal A_\Delta$ of $U^*$, complement of the union
of the hyperplanes of the arrangement.

 From the general theory (cf. \cite{te}), if
  $\Omega_i(\mathcal A_\Delta)$ denotes the space of rational
differential forms of degree $i$ on $\mathcal A_\Delta$ one has 
   the {\it formality},  that is the fact that the $\mathbb
 Z$
subalgebra of differential forms on $\mathcal A_\Delta$ generated by the
linear forms ${1\over 2\pi i}d\log \alpha, \alpha\in\Delta$ is isomorphic (via De Rham theory) to the integral cohomology of   $\mathcal A_\Delta$ \cite{te}.

Formality implies in particular that $\Omega_r(\mathcal
A_\Delta)=H^r\oplus d \Omega_{r-1}(\mathcal A_\Delta),$ for top degree
forms.
$H^r\equiv H^r(\mathcal A_\Delta,\mathbb C)$ is the $\mathbb C$ span of the top degree forms
$\omega_\sigma:=d\log\gamma_1\wedge\dots\wedge d\log\gamma_r$ for all
bases
$\sigma:=\{\gamma_1,\dots,\gamma_r\}$ extracted from $\Delta$.  The forms
$\omega_\sigma$ satisfy a set of linear relations generated by the
following ones. Given $r+1$ elements
$\gamma_i\in\Delta$ spanning $U,$ we have:
$$\sum_{i=1}^{r+1}(-1)^i d\log\gamma_1\wedge\dots
\check{d\log\gamma_i}\dots\wedge d\log\gamma_r=0.$$ 

The projection of  $\Omega_r(\mathcal
A_\Delta)$ to $H^r $  induced by the previous decomposition is by the definition the {\it Total residue} $Tres$.

Recall that a {\it
non broken circuit } in $\Delta$ (with respect to the given total
ordering) is an ordered linearly independent subsequence
$\{\alpha_{i_1},\ldots ,\alpha_{i_t}\}$ such that, for each $1\leq
\ell\leq t$, there is no $j<i_{\ell}$ such that  the vectors
$\alpha_j,\alpha_{i_{\ell}},\dots ,\alpha_{i_t}$ are linearly dependent.
In other words $\alpha_{i_\ell}$ is the minimum element of
$\Delta\cap\langle \alpha_{i_\ell},\ldots ,\alpha_{i_t}\rangle$. In
\cite{te} it is proved that the elements
$$({1\over 2\pi i})^r\omega_{\sigma}:=({1\over 2\pi i})^rd\log\gamma_1\wedge\ldots\wedge d\log\gamma_r,$$ where
$\sigma=\{\gamma_1,\ldots,  \gamma_r \}$ runs over all  ordered bases of
$V$ which are non broken circuits, give a linear $\mathbb Z-$basis of the integral cohomology of   $\mathcal A_\Delta$.\subsection{Irreducibles}

Let us now recall some notions from \cite{dp}. Given a subset $S\subset
\Delta$ we  shall denote by $U_S$ the space spanned by $S$.
\begin{definition} Given a subset $S\subset \Delta $, the completion $\overline
S$ of $S$ equals $U_S\cap \Delta$.  $S$ is called complete if
$S=\overline S$. 

A complete subset $S\subset \Delta $   is called reducible if we can find
a proper partition $S=S_1\dot{\cup} S_2$, called a {\bf decomposition}  such
that 
$U_S=U_{S_1}\oplus U_{S_2}$, irreducible otherwise.

\end{definition} Equivalently we say that the space $U_S$ is reducible. 
Notice that, in the reducible case, $S=S_1\dot{\cup} S_2$,  also $S_1$ and
$S_2$ are complete.

 From this definition it is easy to see \cite{dp}:

\begin{lemma} Given complete sets $A\subset S$ and a decomposition
$S=S_1\dot{\cup} S_2$ of $S$ we have that, if $A=(A\cap S_1)\dot{\cup}(A\cap 
S_2)$ is proper, it is a decomposition of $A$.  Let  $S\subset \Delta$ be  complete. 
Then there is a sequence (unique up to reordering) $S_1,\ldots ,S_m$ of
irreducible subsets in $S$ such that
\begin{itemize}
\item $S=S_1\cup\cdots\cup S_m$ as disjoint union. \item
$U_S=U_{S_1}\oplus\cdots\oplus U_{S_m}.$
\end{itemize} The  $S_i$'s are called the irreducible components of $S$
and the decomposition  $S=S_1\cup\cdots\cup S_m$, the irreducible
decomposition of $S$.
\end{lemma}

 We shall denote by $\mathcal I$ the family of all irreducible subsets in
 $\Delta$.
 
 One of the main points of \cite{dp} is the construction of a smooth variety $X_\Delta$ containing $\mathcal A_\Delta$ as open set with complement a divisor with normal crossing and irreducible components indexed by the set  $\mathcal I$ of irreducibles.    $X_\Delta$  has a proper map to $U^*$ extending the identity of $\mathcal A_\Delta$.
 
 \subsection{Nested sets}  In the theory developed in  \cite{dp} we introduced, for a family of irreducibles $S_i$ the notion of being  {\it nested} according to:
 \begin{definition} A subfamily
$\mathbb M\subset \mathcal I$ is called nested if, given any subfamily
$\{S_1,\ldots ,S_m\}\subset  \mathbb M$ with the property that for no
$i\neq j$, $S_i\subset S_j$, then
$S:=S_1\cup\cdots\cup S_m$ is complete and the $S_i$'s are the
irreducible components of $S$.
\end{definition}  The geometric meaning of this notion is in the fact that, in the model  $X_\Delta$ a set of boundary divisors indexed by a family $\mathbb M\subset \mathcal I$ has non empty intersection if and only if $\mathbb M $ is nested.
We also have \cite{dp1}
\begin{lemma} 1) Let   $\mathbb M=\{S_1,\ldots ,S_m\}$ be a nested set.
Then $S:=\cup_{i=1}^mS_i$ is complete. The irreducible components of $S$
are the maximal elements of  $\mathbb M$.

2) Any nested set is the set of irreducible components of the elements of
a flag $A_1\supset A_2\supset \dots\supset A_k$, where each $A_i$ is
complete.
\end{lemma}

\begin{proposition}  1) Let $A_1\supsetneq A_2\dots\supsetneq A_k$, be a
maximal flag of  complete non empty sets. Then $k=r$ and for each $i$,
$A_i$ spans a subspace of codimension $i-1$.

2)  Let $\Delta=S_1\cup\ldots \cup S_t$ be the irreducible decomposition
of $\Delta$.  

\quad i) Then the $S_i$'s are the maximal elements in $\mathcal I$. 

\quad ii) Every maximal nested set contains each of the elements $S_i$,

\quad\quad  $i=1,\ldots ,t$ and is a union of maximal nested sets in the
sets $S_i$.

3)  Let $ \mathbb M$ be a maximal nested set, $A\in
 \mathbb M$ and $B_1,\dots,B_r\in \mathbb M$   maximal among the elements
in $\mathbb M$ properly contained in $A$. 

Then the subspaces
$U_{B_i}$ form a direct sum and
$$ \dim(\oplus_{i=1}^kU_{B_i})+1=\dim U_A.$$

4) A maximal nested set always has $r$ elements.

\end{proposition}

One way of using the previous result is the following. Given a basis
$\sigma:=\{\gamma_1,\dots,\gamma_r\}\subset \Delta$,  one can associate to
$\sigma$ a maximal flag $F(\sigma)$  by setting
$A_i(\sigma):=\Delta\cap\langle\gamma_{i},\dots,\gamma_r \rangle$. 
Clearly the maps from bases to flags  and from flags to
maximal nested sets are both surjective. We thus obtain a surjective map
from bases to maximal nested sets.  In fact this map induces a
bijection between the set of no broken circuit bases and that of proper maximal nested
sets (see below for their definition).

 We define a map $\phi$ from subsets of $\Delta$ to $ \Delta$ by
 associating to each $S\subset\Delta$ its minimum  $\phi(S):=\min(a\in
S)$ with respect to
 the given  ordering.

We give the definition:
\begin{definition} 1) A flag of complete sets $S_i$ is called proper if
the set  $\phi(S_i)\subset \Delta$ is a basis of $V$.

2) A maximal nested set $\mathbb M$ is called proper if
the set  $\phi(\mathbb M)\subset \Delta$ is a basis of $V$.
\end{definition}    

The main combinatorial result of \cite{dp1} is that:
\begin{theorem}
We have canonical bijective correspondences between:

1) Proper flags of complete sets.

2)  Proper maximal nested sets.

3) Non broken circuit bases.
\end{theorem}
The bijection is given as follows:

Given a basis
$\sigma=\{\gamma_1,\dots,\gamma_r\}$,  we  associate to $\sigma$ the flag
$A_i=\Delta\cap\langle\gamma_i,\dots,\gamma_r\rangle$.

Then  the
maximal nested set   is the decomposition of the previous flag.

\section{A basis for homology}

    Let us denote by
$\mathcal C$ the set of
   non broken circuit  bases of $V$, by  $\mathcal M$ denote the set of
proper maximal nested set.

Let us now fix a basis $\sigma\subset\Delta$. Write
$\sigma=\{\gamma_1,\ldots ,\gamma_r\}$ and consider the $r$-form
$$\omega_{\sigma}:=d\log\gamma_1\wedge\ldots\wedge d\log\gamma_r.$$ This
is a holomorphic form on the open set $\mathcal A_{\Delta}$ of $U^*$
which is the complement of the arrangement formed by  the hyperplanes
whose equation is in $\Delta$. In particular if
$\mathbb M\in \mathcal M$, we shall set $\omega_{\mathbb
M}:=\omega_{\phi(\mathbb M)}$.

Also if $\mathbb M\in \mathcal M$, we can define a homology class in
$H_r(\mathcal A_{\Delta},\mathbb
Z)$ as follows. 
Identify $U^*$ with $\mathbb A^r$ using the coordinates $\phi(S),\ S\in
\mathbb M$. Consider another complex affine space $\mathbb A^r$ with
coordinates $z_S$, $S\in \mathbb M$. In $\mathbb A^r$ take the small
torus $T$ of equation
$|z_S|=\varepsilon$ for each $S\in \mathbb M$.  Define a map
\begin{equation}f:\mathbb A^r\to U^*,\ \ {\rm by} \ \ \phi(S):=\prod_{S'\supset
S}z_{S'}.\end{equation} In \cite{dp} we have proved that this map  lifts, in a
neighborhood of 0, to a local system of coordinates of the model
$X_{\Delta}$  .  To be precise for a vector $\alpha\in\Delta$, set
$B=p_{\mathbb M}(\alpha)$. In the coordinates $z_S$, we have that
\begin{equation}\label{esp}\alpha=\sum_{B'\subset
B}a_{B'}\prod_{S\supseteq B'}z_S=\prod_{S\supseteq
B}z_S(a_B+\sum_{B'\subset B}a_{B'}\prod_{B\supsetneq S\supseteq
B'}z_S)\end{equation} with $a_{B'}\in\mathbb C$ and $a_B\neq 0$.   Set 
$f_{M,\alpha}(z_S):= a_B+\sum_{B'\subset B}a_{B'}\prod_{B\supsetneq
S\supseteq B'}z_S $  and $A_{\mathbb M}$ be the complement in the affine
space $\mathbb A^r$ of coordinates $z_S$ of the hypersurfaces of
equations $f_{M,\alpha}(z_S)=0$. The main point
 is  that 
$A_{\mathbb M}$ is an open set of $X_\Delta$. The point 0 in  $A_{\mathbb
M}$  is the {\it point at infinity} $P_{\mathbb M}$. The open set
$\mathcal A_\Delta$ is contained in    $A_{\mathbb M}$  as the complement
of   the  divisor with normal crossings given by the equations $z_S=0$.
 From this one sees immediately  that
 if $\varepsilon$ is sufficiently small, $f$ maps $T$ homeomorphically
into $\mathcal A_{\Delta}$. Let us give to $T$ the obvious orientation
coming from the total ordering of $\mathbb M$, so that
$H_r(T,\mathbb Z)$ is identified with $\mathbb Z$ and set
$c_{\mathbb M}=f_*(1)\in H_r(\mathcal A_{\Delta},\mathbb Z)$.

Given the class  $c_{\mathbb M}$ and an $r-$dimensional differential form
$\psi$ we can compute $\int_{c_{\mathbb M}}\psi$. Denoting by $P_{\mathbb
M}$ the point at infinity corresponding to 0 in the previously
constructed coordinates $z_i:=z_{S_i}$ we shall say:
\begin{definition} The integral  ${1\over  (2\pi i)^r}\int_{c_{\mathbb
M}}\psi$ is called the {\bf residue} of $\psi$ at the point at infinity
$P_{\mathbb M}$. We will also denote it by $res_{\mathbb M}(\psi)$.
\end{definition}

Notice that the rational forms regular in $\mathcal A_\Delta$,  in a neighborhood of the point 
$P_{\mathbb M}$ and in the coordinates $z_i$, have the form
$\psi=f(z_1,\dots,z_r)dz_1\wedge\dots\wedge dz_r$ with $f(z_1,\dots,z_r)$
a Laurent series which can be explicitly computed  (this  is the consequence of the fact that the model $X_\Delta$ has normal crossings). One then gets that the residue
$res_{\mathbb M}(\psi)$  equals the coefficient of
$(z_1\dots z_r)^{-1}$,  in this series.\smallskip

By abuse of notations, since we have canonical coordinates $z_S$ we shall also speak of residue of a function and write 
$res_{\mathbb M}(f(z_1,\dots,z_r))$

 The
main Theorem of \cite{dp1} is.

\begin{theorem} The set of elements $c_{\mathbb M}$, $\mathbb M\in
\mathcal M$ is the basis of $H_r(\mathcal A_{\Delta}, \mathbb
 Z),$ dual, under the residue
pairing,  to the basis given by the forms
$\omega_{\phi(\mathbb M)}$:  the forms  associated to  the no broken
circuit bases relative to the given ordering.
\end{theorem}
 
 We have seen thus in \cite{dp1} that:
 
1) The formulas found give us an explicit formula for the projection
$\pi$ of $\Omega_r(\mathcal A_{\Delta})=H^r\oplus d
\Omega_{r-1}(\mathcal A_\Delta)$ to $H^r$ with kernel $d
\Omega_{r-1}(\mathcal A_\Delta)$. We have:
\begin{equation}\pi(\psi)=\sum_{\mathbb M\in \mathcal M
}res_{\mathbb M}(\psi)\omega_{\mathbb M}. \end{equation} 

2)  Using the projection $\pi$ any linear map on $H^r$,   can be thought of as a linear map on $\Omega_r(\mathcal A_{\Delta})$
vanishing on $ d
\Omega_{r-1}(\mathcal A_\Delta)$. Our geometric description of homology
allows us to describe any such map as integration on a cycle, linear
combination of the cycles $c_{\mathbb M}$.   
 
\section{The   residue formulas} In this section $V$ is a real
$r-$dimensional vector space and
 $U:=V\otimes_{\mathbb R}\mathbb C$,
$\Delta=\{\alpha_1,\dots,\alpha_n\}\subset V$.  
We fix an orientation for $V$.

 We now further restrict to the case in which  there exists a linear function
on $V$ which is positive on $\Delta$,  i.e. that  all the elements in
$\Delta$ are on the same side of some hyperplane.  The cone   $C$ spanned by the vectors in $\Delta$ is acute and we can decompose it into chambers using the hyperplanes generated by vectors of $\Delta$.  
 
For each  basis $\tau\subset \Delta$,  set $C(\tau)=\{x\in V\,|
x=\sum_{\alpha\in\tau}a_{\alpha}\alpha, a_{\alpha}>0\}$.  Set for
simplicity, for a proper maximal nested set ${\mathbb M}$,
$C({\mathbb M}):=C(\phi({\mathbb M}))$.

 The  final result of \cite{dp1} gives the formula to count the number of integer points  $N_a$  for a polytope $\Pi_a$ in term of residues as follows, choose a big chamber $C$ so that $a\in\overline C$ (the closure of $C$) then: 

\begin{equation} N_a=\sum_{{\mathbb M}\in \mathcal M |C\subset C( {\mathbb
M})}res_{\mathbb M}({e^{<a,x>}\over \prod_{k=1}^N(1-e^{<\alpha_k,x>})}) . \end{equation}

In a similar spirit there is a simpler formula which computes the volume $V_a$ of the polytope as:
\begin{equation}V_a= \sum_{{\mathbb M}\in \mathcal M |C\subset C( {\mathbb
M})}res_{\mathbb M}({e^{<a,x>}\over \prod_{k=1}^N <\alpha_k,x>}) . \end{equation}

The algorithm to compute it takes the following steps:

1. Order the set $\Delta$ and determine the proper maximal nested sets.

2.  If $a\in \overline C$ has been fixed, determine only the proper maximal nested sets such that $C\subset C( {\mathbb
M})$  (often very few out of all the proper nested sets).

3. Prepare for each proper maximal nested set $\mathbb M$ with  $C\subset C( {\mathbb
M})$  the change of new coordinates $z_S $  as in formulas (3), (4) to substitute in the function $ {e^{<a,x>}\over \prod_{k=1}^N(1-e^{<\alpha_k,x>})} $.

4. In the new coordinates $z_S$ each term $1-e^{<\alpha_k,x>}$ equals a product of the variables $z_S$ times an invertible power series in these variables. Hence one can develop enough terms of the function $ {e^{<a,x>}\over \prod_{k=1}^N(1-e^{<\alpha_k,x>})} $ so to be able to compute the residue.

\bigskip
\begin{remark}\label{poli} i)  Each term of formulas (5),(6) is clearly a polynomial in the variable $a$ of which one can estimate the degree.

ii)   If $a$ is a regular vector the condition $C\subset C( {\mathbb
M})$ is equivalent to $a\in   C( {\mathbb
M})$. Otherwise we have in general more than one choice for the chamber $C$ and the formulas are not unique.\end{remark}

Thus these formulas determine functions on the cone $C(\Delta)$  generated by $\Delta$  which are {\it locally polynomials}. More precisely they are polynomials on the strata of an equivalence relation.

Set in fact  $S(v):=\{ \mathbb M \in\mathcal M \ | \, v\in C(\mathbb M)\}.$  Define $v\cong_{\mathcal M} w,\ \iff S(v)=S(w)$ and, if $\mathbb M\in S(v)$, then $v$ and $w$ belong to the same relatively open face of $C(\mathbb M)$. Clearly the strata of this equivalence relation are convex polyhedral cones which decompose $C(\Delta)$. On each of these strata the formulas (5), (6) take polynomial values.  In the appendix we discuss this phenomenon.
\medskip

\section{Network flows arrangements.}
\subsection{Graph arrangements}

The magic arrangement can be seen as a special case of the following general setting.
Let $\Gamma:=(V,L)$  be an oriented graph, i.e. we assume that each $a\in L$ has an initial vertex
$i(a)$ and a final vertex $f(a)$ we also assume that there are no simple loops i.e. edges with initial
and final vertex  equal and that two vertices are joined by at most one edge.

Denote by $l,v$ the number of edges and vertices respectively and by $b_1,b_0$ the two Betti numbers of the graph.
Of course $l-v=b_1-b_0$.

\begin{remark} For a connected graph the number of independent loops is by definition the dimension of  its first homology group, i.e. $l-v+1$.\end{remark}

Often, taking just a subset of the edges we will speak of the graph they generate, meaning that the vertices are exactly all the vertices of the given set of edges.
\smallskip

Fix a basis element $e_v$ for every vertex and consider the set of vectors  $\Delta_\Gamma:=\{x_a:=e_{f(a)}-e_{i(a)}\}$
as $a\in L$.  

\begin{lemma}\label{dimensio} The vectors $x_a$ span a space $V_\Gamma$ of dimension  $v-b_0$.
\end{lemma}
\proof Clearly the spaces spanned by vectors in different connected components form a direct sum so
it suffices to prove the formula when $\Gamma$ is connected.  These vectors span a space $U$   
contained in the subspace generated by the vectors $e_v$ where the sum of the coordinates is 0. We
claim that $U$ coincides with this subspace, in fact choose a vector $e_v$ and add it to $U$ then
by connectedness each $e_w\in U+  Re_v$  hence the claim.\qed
\smallskip

\begin{remark} If $\Gamma$ is connected, $v-1$ edges $a_i$ are such that the vectors $x_{a_i}$ are a basis of $V_\Gamma$ if and only if these edges span a maximal tree.\end{remark}

If we have given a total order to the edges it makes thus sense to ask wether a basis $x_{a_1},\dots,x_{a_{v-1}}$ or a maximal tree $a_1,\dots,a_{v-1}$ is no broken circuit.  

This means that  each $a_i$ is minimal in the complete graph generated by the vertices of  $a_i,\dots,a_{v-1}$.\smallskip

The vectors $\Delta_\Gamma$ define thus a hyperplane arrangement in the $v-b_0$ dimensional space $V_\Gamma^*$ which we shall call a graph arrangement.    We shall now investigate the irreducible subsets in $\Delta_\Gamma$ and the corresponding nested sets.

In case that  the orientation of the   graph is such that  the vectors $x_a$  span an acute cone, i.e. they are all on the same side of a hyperplane which does not contain any of them,  we will speak of  a {\it network} and  we can define the corresponding network polytopes (we simplify re. \cite{BL}).  Let us recall some simple facts about this notion.
\begin{proposition}  A way to obtain a network is by fixing a total order on the set of vertices and orient the edges according to the given order.

An oriented graph is a network if and only if it does not contain oriented loops.

In every network we can totally order the vertices in a way compatible with the orientation of the edges.\end{proposition}
\proof
If we have a total order on the vertices and consider a vector $\alpha$ with strictly increasing coordinates with respect to this total order, we have that the scalar product of $\alpha$ with each $x_a$ is strictly positive.

An oriented loop $a_1,\dots,a_k$ gives vectors $x_{a_1},\dots,x_{a_k}$ with $x_{a_1}+\dots+x_{a_k}=0$ so it cannot be a network.

Conversely assume there are no oriented loops.  Take a maximal oriented chain  $a_1,\dots,a_k$, this is not a loop and necessarily $a_1$ is a source otherwise we could increase the oriented chain.   We take this source to be the smallest vertex, remove it and all the edges coming from it and then start again on the remaining graph by recursion. 
\qed

Remark that $\Delta_\Gamma$ is in canonical bijection with the set of edges of $\Gamma$. Thus subsets of $\Delta_\Gamma$ correspond to subgraphs of $\Gamma$ (with no isolated vertices). Given $A\subset \Delta_\Gamma$ we shall denote by $\Gamma_{A}$ the corresponding graph.

Now recall that a subgraph is called complete if whenever it contains two vertices it also contains all edges between them. On the other hand a subset $A \subset \Delta_\Gamma$ is complete in the sense of arrangements if and only if $\langle A\rangle\cap \Delta_\Gamma=A$. If $A$ is complete in this sense, we shall say that the corresponding subgraph $\Gamma_{A}$  is A-complete. 

\begin{proposition}\label{complete} A connected subgraph of  $\Gamma$ is $A$-complete if and only if it is complete.

A subgraph is $A$-complete if and only if all its connected components are $A$-complete.\end{proposition}
\proof The fact that a complete subgraph is also A-complete is clear.

If  a subgraph $\Lambda$ is not A-complete we have an edge $a\notin\Lambda$ which is dependent of the edges in $\Lambda$.  We know that the
dimension of the corresponding span equals the number of vertices in $\Lambda$ minus the number of its connected components. This dimension can remain the same if and only if, adding $a$, we do not add any vertices nor do we decrease
the number of connected components.    If $\Lambda$ is connected, this means that the vertices of $a$ are in $\Lambda$.  Hence $\Lambda$ is not complete.  

If the graph is not connected the condition is not only that the vertices of $a$ are in $\Lambda$ but also that they belong to the same connected component, hence the claim.

\qed

The previous proof has a simple but important consequence:
\begin{corollary}
A  graph arrangement is unimodular.
\end{corollary}

\begin{corollary}\label{maximal} Given a connected graph $\Gamma$, a proper subgraph $\Lambda$ is maximal $A$-complete
if and only if either  $\Lambda$ is a connected subgraph obtained from $\Gamma$  deleting one vertex and all the edges from it, or it is a graph with two connected components obtained from $\Gamma$ deleting a set of edges each of which joins the two components.\end{corollary}
\proof If we remove one vertex and all the edges from it,  and the resulting graph $\Lambda$, with edges $C$,   is still connected, it follows that the corresponding subspace $<C>$ has codimension 1. Since $\Lambda$ is clearly complete it is also maximal.   Similarly in the second case where the graph we obtain is complete the number of vertices is unchanged but it has two connected components, so it gives again a codimension 1 subspace.

Conversely if $\Lambda$ is maximal $A$-complete
with $w$ vertices and $b$ connected components we must have  $v-2=w-b$  so, either $b=1$ and $w=v-1$ or   $w=v$ and $b=2$. It is now easy to see that we must be in one of the two preceding cases ( from the description of  $A$-complete subgraphs and Lemma \ref{dimensio}).\qed

 A set of edges so
that the remaining graph has two connected components and all the deleted edges join the two components will be called a {\it simple disconnecting set}.  To find such a set is equivalent as to decompose the set of vertices into two disjoint subsets  $V_1,V_2$, so that each of the two complete subgraphs of vertices $V_1,V_2$ are connected.
 \smallskip

Recall that, given two graphs $\Gamma_1,\Gamma_2$ with a preferred vertex $w_1,w_2$ in each, the {\it wedge $\Gamma_1\vee \Gamma_2$ } of the two graphs is given by forming their disjoint union and then identifying the two vertices.  Clearly if $v_1,v_2,v$ resp. $b_1,b_2,b$ denote the the number of vertices resp. of connected components of $\Gamma_1,\Gamma_2,\Gamma_1\vee \Gamma_2$ we have $v=v_1+v_2-1,\ b=b_1+b_2-1$  hence  $v-b=v_1-b_1+v_2-b_2$

We shall now say that an $A$-complete subgraph is irreducible if the corresponding subset of $\Delta_\Gamma$ is irreducible. The previous formulas shows that the decomposition of a graph as wedge or into connected components, implies a decomposition of the corresponding hyperplane arrangement so that in order to be irreducible a graph must be connected and cannot be expressed as the wedge of two smaller subgraphs.  The following Proposition  shows that these conditions are also sufficient.

\begin{proposition}\label{irredu}  A connected graph $\Gamma$  is irreducible if and only if it is not a wedge of two graphs. This also means that there is no vertex which disconnects the graph.
\end{proposition}
\proof
We have already remarked  that the decomposition of a graph as wedge implies a decomposition of the corresponding hyperplane arrangement.

On the other hand, let $\Gamma$ be  connected  and suppose that $\Delta_{\Gamma}$ has a non trivial decomposition $\Delta_{\Gamma}=A\cup B$ (in the sense of hyperplane arrangements). If $A'$ is the set of edges of a connected component of $\Gamma_A$, then $\Delta_{\Gamma}=A'\cup ((A-A')\cup B)$ is also a decomposition so we can assume that   $\Gamma_A$ is connected. 

Denote by $V_A$ (resp. $V_B$) the set of vertices of  $\Gamma_A$  (resp. $\Gamma_B.$)
We must have that $V_A\cap V_B$ is not empty since $\Gamma$ is connected. 

The fact that $\Delta_{\Gamma}=A\cup B$ is a decomposition implies $\langle\Delta_{\Gamma}\rangle=\langle A\rangle\oplus\langle B\rangle$ so from Lemma \ref{dimensio} we deduce   $v-1=|V_A|-b_A+|V_B|-b_B$  (with $b_A,b_B$ the number of connected components of the two graphs with edges $A,B$ respectively).   Since $v=|V_A|+|V_B|-|V_A\cap V_B|$,  we have $1+ |V_A\cap V_B|=  b_A+ b_B$.  

We are assuming that $b_A=1$, so we get  $ |V_A\cap V_B|=    b_B$. Since $\Gamma $ is connected each  connected component  of    $\Gamma_B$ must contain at least one of the vertices in  $V_A\cap V_B$.  The equality  $ |V_A\cap V_B|=    b_B$ implies then that each  connected component  of    $\Gamma_B$   contains exactly  one of the vertices in  $V_A\cap V_B$. Thus $\Gamma $ is obtained from $\Gamma_A$ by attaching, via a wedge operation each  connected component  of    $\Gamma_B$ on different vertices.
\qed

\begin {remarks} 1) Notice that in fact, the first case of \ref{maximal} could be considered as a degenerate case of the second, with $\Gamma_A$ reduced to a single vertex.
\smallskip

2) In general a  complete decomposition will thus present a connected graph as an iterated wedge of irreducible graphs.\end{remarks}
 
With the previous analysis it is easy to give an algorithm which allows to describe all proper maximal nested sets.

The algorithm is recursive and based on the idea of building a proper maximal flag of complete sets.

Step 1. We choose a total order of the edges.
 
 Step 2. We decompose the graph into irreducibles.

Step 3.  We proceed by recursion separately on each irreducible where we have the induced total order.

We assume thus we have chosen one irreducible.

Step 4. We build all the proper maximal complete sets which do not contain the minimal edge. These are of two types. 

i) The two sets obtained by removing the edges out of one of the vertices of the minimal edge (by Proposition \ref{irredu} this operation produces a connected graph).

ii) Remove the simple disconnecting sets containing the minimal edge.

Step 5.  Keeping the induced order, go back to Step 2 for each of the proper maximal complete sets constructed.
\bigskip

 From the residue formulas it is clear that, in the previous algorithm, given a vector $u$ in order to compute the volume or the number of integer points of the polytope $\Pi_a$, it is only necessary to compute the proper nested sets $\mathbb M$ which satisfy the further condition $u\in C(\mathbb M)$, which we will express by the phrase {\it $\mathbb M$  is adapted to $u$}.    The previous algorithm explain also how to take into account this condition.

In fact let  $\Lambda\subset \Gamma$ be a proper maximal complete set, we can see as follows if this can be the first step to construct an  $\mathbb M$  adapted to $u$. In fact the basis  $\phi(\mathbb M)$ is composed of the minimal element $x_a$ and the basis of the span $<\Lambda>$  corresponding on the part $\mathbb M'$ of the nested set contained in $<\Lambda>$. Thus $u$ can be written uniquely in the form  $\lambda x_a+w$ with $w\in <\Lambda>$.

$\mathbb M$ is adapted to $u$ if and only if $\lambda\geq 0$ and $\mathbb M'$ is adapted to $w$.

This gives a recursive way of proceeding if we can compute $\lambda$. Let us do this in the second case (since the first is a degenerate case).  In the decomposition of the maximal complete subset as $A\cup B$ let us assume that the orientation of the arrow of the minimum edge $a$ points towards $A$ so that $x_a$  as a function on the vertices takes the value 0 on all vertices except its final point in $A$ where it takes the value 1, and its initial point in $B$  where it takes the value - 1.  The vector $u$  is just a function on the vertices with the sum 0.  Let $\lambda$ equal the sum of the values of $u$ on the vertices of $A$, thus    $-\lambda$ equals the sum of the values of $u$ on the vertices of $B.$

We then have that $u-\lambda x_a=w\in <\Lambda>$, so we see that:
\begin{proposition} In the decomposition of the maximal complete subset $\Lambda $ as $A\cup B$ let us assume that the orientation of the arrow of the minimum edge $a$ points towards $A$. Then if $\Lambda$ is the first step of a proper flag adapted to $u$ we must have that the sum of the values of $u$ on the vertices of $A$ is non negative.

\end{proposition}
{\bf Remark} We have seen in 3.6  that the decomposition into chambers can be detected by any choice of ordering and the corresponding n.b.c. bases.  Nevertheless  the number of n.b.c. bases adapted to a given cell depends strongly on the order.  For a given cell it would thus be useful to minimize this number in order to  optimize the algorithms computing formulas  (5),(6). This point needs a further investigation which we have not done.

\subsection{Two examples: $A_n$ and  Magic arrangements} As we have mentioned  before, in the case our graph $\Gamma$ is the complete graph on $n+1$ elements $\{1,\ldots ,n+1\}$,  the arrangement we obtain is the root arrangement of type $A_n$. If we furthermore order  our vertex set in the obvious way we get that $\Delta_{\Gamma}=\{e_i-e_j, \ 1\leq i<j\leq n+1\}$ is the set of positive roots. This case has been studied extensively (see for example \cite{dp1}). Our previous analysis allows us to recover immediately a number of know facts. Given any set $S$ of vertices with at least two elements, the complete subgraph with vertex set $S$ is clearly a complete graph and hence irreducible. It follows that   irreducible subsets of $\Delta_{\Gamma}$ are in bijection with subsets of $\{1,\ldots ,n\}$ containing at least 2 elements.

Under this correspondence,  a sequence $S_1,\ldots , S_t$ of subsets of $\{1,\ldots ,n+1\}$ containing at least 2 elements is nested if and only
for any $1\leq i,j\leq t$, either $S_i\cap S_j=\emptyset$ or $S_i$ and $S_j$ are one contained in the other. 

Let us now fix the following total order on $\Delta_{\Gamma}$. We set $e_i-e_j\le  e_h-e_k$ if $k-h<j-i$ and if $k-h=j-i$, if $i\le h$. 
A proper maximal nested set  $\mathbb M$ is then encoded by a sequence 
of $n$ subsets each having at least two elements,  with the property
that, taking the minimum and maximum for each set, these pairs are all
distinct.

Correcting an imprecision in \cite{dp1} pag.5 let us explain  how to inductively define a bijection between
proper maximal nested sets   and  permutations  of $1,\dots,n$ fixing
$n$.  To see this consider a maximal nested set  $\mathbb M$ as a sequence $\{S_1,\ldots,
S_{n}\}$ of subsets of $\{1,\ldots ,n+1\}$ with the above properties.  
We can assume that 
 $S_1=(1,2,\ldots,         n+1)$. Using Corollary \ref{maximal}  we see that $\mathbb M':=\mathbb M-\{S_1\}$
has either one or two maximal elements. 
  If $S_2$ is the unique maximal element and $1\notin S_2$, by
 induction we get  a  permutation $p(\mathbb M')$ of $2,\ldots ,n$ fixing $n$. We
then  set $p(\mathbb M)$ equal to the permutation which fixes $1$ and is
equal to $p(\mathbb M')$ on  $2,\ldots ,n$. 
  If $S_2$ is the unique maximal element and $n\notin S_2$, 
 we get, by induction, a  permutation $p(\mathbb M')$ of $1,\ldots ,n-1$ fixing $n-1$. We then  set
$p(\mathbb M)$ equal to the permutation which fixes $n$ and is equal to
$\tau p(\mathbb M')$ on $S_2=\{1,\ldots ,n-1\}$,  $\tau$ being  the
permutation which reverses the order in $S_2$.
   If $S_2$ and $S_3$ are the two maximal elements  so that  $\{1,\ldots
,n\}$ is their disjoint union,  and $1\in S_2=\{1=i_1<\cdots <i_h \}$,
$n\in S_3=\{i_{h+1} <\cdots <i_n=n\}$ then by induction we get two permutations $p_2$ and $p_3$ of
$S_2$ and $S_3$ respectively.   A permutation $\sigma$ of a subset $S=\{i_1<\cdots <i_h\}\subset \{1,\ldots ,n\}$ induces a bijection $\sigma'$ between $\{1,\ldots ,h\}$ and $S$ defined by $\sigma'(t)=\sigma(i_t)$ for $1\leq t\leq h$. We  then set 
 $p(\mathbb M)(t)$ equal to $p_3'(t-|S_2|+1)$ if $t>|S_2|$  and equal to 
$(\tau p_2)'(t)$ otherwise,  $\tau$ being  the permutation which reverses
the order in $S_2$.   Remark that the two sets  $S_2,S_3$  are determined by $p(\mathbb M) $ by writing it as a word and collecting all the entries appearing before and including 1 and all the entries after.   The two permutations are also similarly reconstructed. In particular this shows that there are $(n-1)!$ 
proper maximal nested sets, which can be recursively constructed. 

\smallskip 

The second example we want to analyze is the following:

Given  2 positive integers $m,n$ we define the arrangement $M(m,n)$ as follows. We start from the vector space $\mathbb R^{m+n}$ withe basis elements $e_i,\ i=1,\dots,m,\ f_j,\ j=1,\dots, n$ and let $V$ be the hyperplane where the sum of the coordinates is 0. The arrangement is given by the $nm$ vectors $\Delta(m,n):=\{(i|j):=e_i-f_j, \ i=1,\dots,m,\  \ j=1,\dots, n\}$. It is the graph arrangement associated to the full bipartite graph formed of all oriented edges from a set $X$ with $n$ elements to a set $Y$ of $m$ elements.\smallskip

Let us discuss the  notions of irreducible, nested and proper nested in the example $M(m,n)$. We need some definitions,  given two non empty subsets  $A\subset\{1,2,\dots,m\}, B\subset\{1,2,\dots,n\},$ we denote by $A\times B$ the set of vectors $(i|j), i\in A,\ j\in B$ and call it a rectangle.  We say that the rectangle is degenerate if  either $A$ or $B$ consists of just one element (and we will speak of a row or a column respectively).
 
 In particular when $A,B$ have two elements we have a {\it little square}. We define {\it triangle}   a subset of 3 elements of a little square.

 \begin{lemma}  
1)  The 4 elements of a little square  $\{i,j\}\times\{h,k\}$ form a complete set. They span a 3 dimensional space and satisfy the relation $(i|h)+(j|k)=(i|k)+(j|h)$
 
 2) The completion of a triangle is the unique little square in which it is contained.
 
 3) Any rectangle is complete.
 \end{lemma} 
 The proof is clear.

 \begin{theorem}\label{magiccom}
 For a subset $S\subset \Delta(m,n)$ the following conditions are equivalent.
 
 1) $\Gamma_S$ is A-complete.
 
 2)  If a triangle $T$ is contained in $S$ then its associated little square is also contained in $S$.
 
 3) $S=\cup_{i=1}^h A_i\times B_i$ where the $A_i$ are mutually disjoint and also the $B_j$ are mutually disjoint .
 \end{theorem}
 \proof
 Clearly 1) implies 2).  Assume 2), consider a maximal rectangle $A\times B$ contained in $S$, we claim that  $S\subset A\times B\cup C(A)\times C(B)$. Otherwise there is an element $(i|k)\in S$ where either $i\in A, k\notin B$ or  $i\notin A, k\in B$. Let us treat the first case the second is similar.  If $A=\{i\}$ then $A\times (B\cup\{k\})$ is a larger rectangle contained in $S$ a contradiction. Otherwise take $j\in A, j\neq i, h\in B$ we have that  $(i|h),  (j|h), (i|k)$ are in $S$ and form a triangle so by assumption also $(j|k)\in S$ this means that again  $A\times (B\cup\{k\})$ is a larger rectangle contained in $S$ a contradiction.   Now we can observe that  $S\cap 
C(A)\times C(B)$ is also complete and we proceed  by induction. 

3) implies 1) follows from Proposition \ref{complete}.\qed

Theorem \ref{magiccom} now gives  the decomposition of a complete set  into irreducibles.

\begin{corollary}
A non degenerate rectangle is irreducible. Given a complete set of the form $S=\cup_{i=1}^h A_i\times B_i$ where the $A_i$ are mutually disjoint and also the $B_j$ are mutually disjoint its irreducible components are the non degenerate rectangles $A_i\times B_i$ and the single elements of the degenerate rectangles $A_i\times B_i$.
\end{corollary}
Theorem \ref{magiccom} also implies  the structure of the maximal proper complete subsets of $\Delta(h,k)$.
\begin{corollary}  A  maximal proper complete subset $S$ of $\Delta(h,k)$ is of one of the following types:

$A\times B\cup C(A)\times C(B)$ if $A,B$ are proper subsets.

$A\times \{1,\dots,n\}$ where $A$ has $m-1$ elements.

$ \{1,\dots,m\} \times B$ where $B$ has $n-1$ elements.

\end{corollary}

All these considerations allow us to find all  proper flags in the case of the magic arrangement. Of course in order even to speak about proper flags, we have to fix a total ordering among the pairs $(i,j)$.  Let us use as order the lexicographic order  so that $(1|1)$ is the minimum element. It follows that if $S$ is a proper maximal complete subset, in order to be the beginning of a proper flag one needs that  $(1|1)\notin S$.

It is then clear that, once we have started with such a proper maximal complete subset we can complete the flag to a proper flag by taking a proper flag in $S$ for the induced order.  This gives a recursive formula for the number $b(m,n)$ of proper flags which is also the top Betti number. We have from our discussion the recursive formula for $b(m,n)$:
$$ \!\!\!b(m-1,n)+b(m,n-1)+\sum_{a,c}\binom{m-1}{a-1}\binom{n-1}{c}b(a,c)b(m-a,n-c) .$$

\section{APPENDIX}
\subsection{The big cells}

   The   condition of \ref{poli} is clearly independent of the order chosen so that it makes sense to ask whether  the   stratification discussed in that remark is independent of the order chosen.   In order to prove this we need a few simple combinatorial lemmas   on polytopes whose proof we recall for completeness.

In an $r-$dimensional real vector space $V$, let us choose a finite set of vectors $\Psi:=\{v_i\}$ spanning $V$ and   lying in an  affine hyperplane $\Pi$ of equation $\langle \phi,x\rangle=1$ for some linear    form $\phi$ . 

 The intersection of  the cone $C(\Psi)$ with $\Pi$ is the convex polytope $\Sigma$ envelop of the vectors $v_i$. Each  cone, generated by $k+1$ independent vectors in $\Psi$,  intersects $\Pi$ in  a $k$ dimensional simplex.   Then the configuration of cones is obtained by projecting a configuration of simplices and there is a simple dictionary to express properties of cones in terms of simplices and conversely.
 
It is well known (and in any case will follow from our more precise results)  that $\Sigma$ is the union of the simplices with vertices independent vectors of $\Psi$.    It is natural to define {\it regular} a point in $\Sigma$  which is not contained in any $r-2$ dimensional simplex (or in the corresponding cone).   The connected components   of the set of regular points are called in \cite{BV1} the {\it big cells}. They are the natural loci where the formulas (5),(6) take polynomial values.  Since on the other hand in \ref{poli}  the natural strata are the ones determined via n.b.c. bases it is important to compare the two stratifications. Our main result is in fact that they coincide  (Theorem \ref{coincide}), this gives a rather strong simplification in the algorithms necessary to determine the big cells.
 In order to do this we have to work in a slightly more general setting  and    consider the following stratifications of $\Sigma$.

Let us choose a family $\mathcal I$ of bases extracted from $\Psi$. For each such basis $\mathfrak b$ consider the family $\mathcal F_{\mathfrak b}$ of  faces of the simplex $\sigma_{\mathfrak b}$ generated by $\mathfrak b$ and set $\mathcal F_{\mathcal I}=\cup_{\mathfrak b\in\mathcal I}\mathcal F_{\mathfrak b}$.  Given a point  $v \in \Sigma$ we set $Z(v)=\{f\in \mathcal F_{\mathcal I}|v\in f\}$ and define an equivalence relation $R_{\mathcal I}$  on $\Sigma$ by setting $v$ and $w$ as equivalent (or belonging to the same stratum), if $Z(v)=Z(w)$.

We want to compare the equivalence relations $R_{\mathcal I}$ for various choices of $\mathcal I$. We start with a special case assuming that the set $\Psi$ consists of  $r+1$ vectors $v_0,\dots,v_{r}$.   We set $\mathcal I$ equal to the family of all bases formed by elements in $\Psi$ and chosen $0\leq j\leq r$, $\mathcal I_j$ equal to the family of all bases formed by elements in $\Psi$ and which contain the vector $v_j$.  Then

\begin{lemma}\label{uno}
$\Sigma =\cup_{\mathfrak b\in\mathcal I_j}\sigma_{\mathfrak b}$,
($\Sigma =\cup_{\mathfrak b\in\mathcal I}\sigma_{\mathfrak b}$).
\end{lemma}

\proof
By suitably reordering we can assume that $j=r$. If $v_0,\ldots v_{r-1}$ are not linearly independent, let us  consider their convex envelope $\Sigma'$ which is contained in the hyperplane $V'$ which they span. By induction we can assume that $\Sigma'$ is the union of the $r-2$- dimensional simplices $\sigma_{\mathfrak b'}$, where $\mathfrak b'$ runs  over the bases of $V'$ which can be extracted from $\{v_0,\ldots ,v_{r-1}\}$. In this case $\Sigma$ is a pyramid with basis  $\Sigma'$ and  and vertex $v_{r}$, and  our claim follows. 

Let us now suppose that $\mathfrak b=\{v_0,\ldots v_{r-1}\}$ is a basis of $V$. Take $v\in \Sigma$. Consider the  line joining $v$  with $v_{r}$. It intersects   $\sigma_{\mathfrak b}$ in a segment with ends two points $a,b$ and then $v$ is either in the segment $av_{r}$ or in $bv_{r}$.    $a,b$ are in $r-2$ dimensional faces of  $\sigma_{\mathfrak b}$,  thus  $v$ lies in the convex envelop of an $r-2$ dimensional face $\tau$ of  $\sigma_{\mathfrak b}$ and $v_r$. If $v_r$ is independent of $\tau$ we have thus an $r-1$ dimensional simplex  having $v_{r}$ as a vertex in which $v$ lies.
Otherwise  by induction  $v$ lies in any case in some simplex  having $v_{r}$ as a vertex which is then contained in a larger $r-1$ dimensional simplex.  \qed

\begin{lemma}\label{unip}
Let $\sigma$ be a simplex, $q$ a point in the interior of a face $\tau$ of $\sigma$ and $p$ a point. Assume that the segment $pq$ intersects $\sigma$ only in $q$ then  the convex hull of $\tau$ and $p$ is a simplex and meets  $\sigma$ in $\tau$. \end{lemma} 
\proof  If $p$ does not lie in the affine space spanned by $\sigma$ the statement is obvious. Otherwise we use the  vector notations, we can assume that $\sigma$  is the convex hull of  the basis vectors $e_1, \dots,e_m$ and $\tau$ is the face of vectors with non zero (positive) coordinates for $i\leq k$.

The condition that $\{t q+(1-t)p\ |\ 0\leq t\leq 1\} \cap\sigma=\{q\}$ is equivalent  to the fact that there is a $i$ larger than $k$ such that the $i-$th coordinate of $p$ is  negative.  This  condition  does not depend on the point $q\in \tau$ and shows that $p$ is independent of $\tau$. Our claim follows.
\qed

We then have 

\begin{lemma}\label{uni}
With the same notation as above, the equivalence relations $R_{\mathcal I}$ and $R_{\mathcal I_j}$ coincide. \end{lemma} 
\proof
As before assume $j=r$. If $\mathfrak b=\{v_0,\ldots v_{r-1}\}$ is not a basis, $\mathcal I=\mathcal I_j$ and there is nothing to prove. Otherwise $\mathcal F_{\mathcal I}$ differs from $\mathcal F_{\mathcal I_j}$ only by adding the interior ${\sigma}^{0}_{\mathfrak b}$ of the simplex $\sigma_{\mathfrak b}$.

Thus we only have two show that  if two vectors  $v,w$ are congruent under $R_{\mathcal I_j}$ it is not possible that $v\in{\sigma}^{0}_{\mathfrak b}$ and  $ w\notin {\sigma}^{0}_{\mathfrak b}$.   
The half line  joining $v_{r}$ and $w$ meets for the first time the simplex $\sigma_{\mathfrak b}$ in a point $u$  which is in the interior of some face $\tau$  of $\sigma_{\mathfrak b}$ and $w$ is in the  segment $v_{r},u$.   By the previous lemma the convex hull of $v_{r}$ and $\tau$  is a simplex and meets $\sigma_{\mathfrak b}$  exactly in $\tau$.  
By hypothesis then $v$ also lies in this simplex hence not in the interior of $\sigma_b$ a  contradiction. 

 \qed

  Let us now go back to our set of vectors $\Delta$ and let us fix an ordering of $\Delta$. Since we know that there is a linear form $\phi$ which takes positive values on each element in $\Delta$, by suitably rescaling with positive numbers, we can assume that $\langle\phi , \alpha\rangle=1$ for each $\alpha\in \Delta$.  We now set $\mathcal I$ equal to the family of all  bases which can be  extracted from $\Delta$ and,as before, $\mathcal M$ equal to the family of n.b.c  bases with respect to the chosen ordering. 

Choose an element $\alpha\in\Delta$ and assume that the element $\beta$ is the successor of $\alpha$ in our ordering. Define a new ordering by exchanging $\alpha$ and $\beta$. The following Lemma tells us how the set  $\mathcal M$ changes.
  \begin{lemma}\label{nobr}  A n.b.c. basis  $\sigma:=\alpha_1,\dots,\alpha_n$ for the first order remains n.b.c. for the second unless all the following conditions are satisfied:
\begin{enumerate}[i)]
\item $\alpha=\alpha_i$ appears in $\sigma$.
\item $\beta$ does not appear in  $\sigma$.
\item $\beta$ is dependent on $\alpha=\alpha_i$ and the elements $\alpha_j, j>i$ in  $\sigma$ following $\alpha$.
\end{enumerate}
\noindent In all these conditions hold,  $\sigma':=\alpha_1,\dots,\alpha_{i-1},\beta,\alpha_{i+1},\dots,\alpha_n$ is a n.b.c. basis for the second order. All   n.b.c bases for  the second order are obtained in this way.\end{lemma}
\proof The proof is immediate and left to the reader.\qed
\begin{theorem}\label{coincide}
The equivalence relations $R_{\mathcal I}$ and $R_{\mathcal M}$ coincide.
\end{theorem}
\proof
As before choose an element $\alpha\in\Delta$ and call $\beta$  the successor of $\alpha$ in our ordering. Define a new ordering by exchanging $\alpha$ and $\beta$ and denote by  $\mathcal M'$ the family of n.b.c  bases with respect to the new  ordering. We claim that the equivalence relations $R_{\mathcal M}$ and $R_{\mathcal M'}$ coincide.

Since every basis extracted from $\Delta$ is a n.b.c. basis for a suitable ordering and we can pass from one ordering to another by a sequence of elementary moves consisting of exchanging an element with its successor, this will prove our Theorem.

Set $\overline{\mathcal M}=\mathcal M\cup\mathcal M'$. Take a basis $\mathfrak b\in \overline{\mathcal M}-\mathcal M$. By Lemma \ref{nobr}
$\mathfrak b=\{\gamma_1,\ldots,\gamma_{k-1} ,\beta,\gamma_{k+1}\ldots ,\gamma_r\}$ with $\alpha,\beta,\gamma_{k+1}\ldots ,\gamma_r$ linearly dependent. Consider the set of vectors $\mathfrak b\cup\{\alpha\}$. To this set we can apply Lemma \ref{uni} and deduce that  the equivalence relation induced by the family of all bases extracted from $\mathfrak b\cup\{\alpha\}$ coincides with the equivalence relation induced by subfamily of all bases containing $\alpha$. These are easily seen to lie all in $\mathcal M$. We deduce that  $R_{\mathcal M}$ and $R_{\overline{\mathcal M}}$ coincide. By symmetry  $R_{\mathcal M'}$ and $R_{\overline{\mathcal M}}$ coincide too, hence our claim.\qed

Notice that  given $v,w\in C(\Delta)$,  we have $v\cong_{\mathcal M}w$ if and only if either $v=w=0$ or $v/\langle\phi,v\rangle\simeq_{R_\mathcal M}w/\langle\phi,v\rangle$. By this remark and our Theorem  it is immediate to see that:
\begin{corollary} If we remove from $C(\Delta)$ the strata which are not of maximal dimension, the resulting connected components are just the big chambers as defined in \cite{BV1}.\end{corollary}

\end{document}